 \documentclass[final,1p,times]{elsarticle}
\usepackage{listings}
\usepackage{color}
\usepackage[]{multicol}
\usepackage{amsthm}
\usepackage{caption}
\usepackage{here}
\usepackage{bm}
\usepackage{mathrsfs}
\usepackage{subcaption}
\usepackage{amsmath, amssymb}
\usepackage{ascmac}
\usepackage{url}
\usepackage{comment}
\usepackage{nccmath}
\usepackage {booktabs}

\newtheorem{theorem}{Theorem}

\newtheorem{lemma}{Lemma}
 
\newtheorem{corollary}{Corollary}
\newtheorem{proposition}{Proposition}
\newtheorem{remark}{Remark}
\def\vec#1{\mbox{\boldmath{$#1$}}}
\journal{}
\begin{document}

\begin{frontmatter}



\title{Exact Distribution of the Noncentral Complex \\Roy's Largest Root Statistic via Pieri's Formula
}


\author[label1]{Koki Shimizu}
\author[label1]{Hiroki Hashiguchi}

\address[label1]{Tokyo University of Science, 1-3 Kagurazaka, Shinjuku-ku, Tokyo, 162-8601, Japan}
\cortext[mycorrespondingauthor]{Corresponding author. Email address: \url{k-shimizu@rs.tus.ac.jp}~(K. Shimizu).}

\begin{abstract}
In this study, we derive the exact distribution and moment of the noncentral complex Roy's largest root statistic, expressed as a product of complex zonal polynomials.
We show that the linearization coefficients arising from the product of complex zonal polynomials in the distribution of Roy's test under a specific alternative hypothesis can be explicitly computed using Pieri's formula, a well-known result in combinatorics.
These results were then applied to compute the power of tests in the complex multivariate analysis of variance (MANOVA).
\end{abstract}




\end{frontmatter}


\section{Introduction}
The matrix variate beta distributions of type I and II are generalizations of the univariate beta and $F$ distributions, respectively.
Matrices following these distributions are referred to as Beta and $F$ matrices, respectively, and their eigenvalues are used for test statistics in multivariate hypothesis testing, such as Roy's largest root test, Pillai-Bartlett trace, Lawley--Hotelling trace, and Wilks'lambda.
Matrix variate distributions have been developed primarily in the real-valued setting, and recent studies have extended them to complex case with applications in multivariate analysis and signal processing.

In the noncentral case, the eigenvalue distributions of the Beta and $F$ matrices play a central role in obtaining the power of various hypothesis tests. 
However, numerical computation under a generalized alternative hypothesis is more difficult than in the central case.
Therefore, exact and approximate distributions have been derived under specific alternative hypotheses.
For the complex case, Gia and Phong \cite{thu2022distribution} derived the exact distribution of the noncentral Wilks statistic. 
Phong \cite{thanh2023approximation} provided a saddlepoint approximation of this distribution.  
Roy's test can be considered more powerful than Wilks' Lambda test when the noncentrality matrix has only one nonzero eigenvalue \cite{butler2020mimo}.
Thus, in Roy's test, the primary interest is in assessing the power under a linear alternative, which is the noncentral matrix that has only one eigenvalue. 
Dharmawansa et al. \cite{dharmawansa2019roy} provided an approximate distribution of the noncentral complex Roy's largest root statistic under a linear alternative.
Dharmawansa et al. \cite{dharmawansa2024detection} derived the exact distribution of the noncentral Roy's largest root statistic under a linear alternative using a determinant representation and applied it to signal detection problems.
However, the approximate distribution of Roy's statistic is not sufficiently accurate in high dimensions, and the computation of the exact distribution suffers from numerical instability during the evaluation of determinants.

In this study, we extended the exact distribution result of the noncentral Roy's largest root statistic in the real case as provided by Shimizu and Hashiguchi \cite{shimizu2022numerical} to the complex case.
In Section~\ref{sec:02}, we discuss the coefficients arising when expressing the product of complex zonal polynomials as linear combinations of complex zonal polynomials.
Additionally, we introduce classical combinatorial tools, including the Littlewood--Richardson rule and Pieri's formula, which are useful for computing products of complex zonal polynomials. 
In Section~\ref{sec:03}, the exact distribution and moments of the noncentral Roy's largest root statistic are derived using the products of complex zonal polynomials.
Section~\ref{sec:04} presents the numerical computation of the distribution using Pieri's formula.

\section{Linearization coefficients of product of complex zonal polynomials}
\label{sec:02}
In this section, we discuss the linearization coefficients that appear when expressing the product of complex zonal polynomials as a linear combination.
For a positive integer $k$, let $\kappa=(\kappa_1,\kappa_2,\dots,\kappa_l)$ denote a partition of $k$ with $\kappa_1\geq\cdots\geq \kappa_l> 0$ and $\kappa_1+\cdots +\kappa_l=k$, where $l=l(\kappa)$ is the length of partitions.
It is convenient to consider $\kappa$ as having any number of additional zeros, $\kappa=(\kappa_1,\kappa_2,\dots,\kappa_l, 0,\dots, 0)$.
The set of all partitions with lengths less than or equal to $m$ is denoted by $P^k_{m}=\{\ \kappa=(\kappa_1,\dots,\kappa_m)\mid \kappa_1+\dots+\kappa_m=k, \kappa_1\geq \kappa_2\geq\cdots\geq  \kappa_m \geq 0 \}$.
  
For $\kappa \in P^k_m$ and $\tau \in P^t_m$, we express the product of Jack polynomials $C^{\beta}(X)$, which are well-known symmetric and homogeneous polynomials associated with an $m \times m$ Hermitian matrix $X$, as a linear combination of Jack polynomials:
\begin{align}
\label{prod-zonal}
C^{\beta}_\kappa(X) \cdot C^{\beta}_\tau(X) = \sum_{\delta \in P^{k+t}_m} g^\delta_{\kappa, \tau}(\beta)\, C^{\beta}_\delta(X),
\end{align}
where $g^\delta_{\kappa, \tau}(\beta)$ is a constant determined by the partitions $\kappa$ and $\tau$.  
For further details on Jack polynomials, see Stanley \cite{stanley1989some} and Macdonald \cite{macdonald1998symmetric}.
If $\beta=1$ and $\beta=2$, $C^{\beta}_\kappa(X)$ is the real and complex zonal polynomials, respectively. 
In this study, we only consider $\beta = 2$, and simply write $C_\kappa(X):=C^{(2)}_\kappa(X)$ for the complex zonal polynomials. 
Although the notation $C_\kappa(X)$ is typically used for the real case, we adopt this notation for the complex case because the real case is not considered.
For the complex zonal polynomials, 
Khatri \cite{khatri1970moments} provided the explicit expression for $C_\kappa(X)$ with $X=I_m$ by
\begin{align*}
C_{\kappa}(I_{m})=\chi^2_\kappa
\frac{\Gamma_m(m, \kappa) }{ \Gamma_m(m)k!},
\end{align*}
 where
$ \Gamma_m(m, \kappa) = \pi^{m (m-1)/2} \prod_{i=1}^{m} \Gamma(m + \kappa_i - i + 1)$, 
$\Gamma_m(m) = \Gamma_m(m, 0)$
and $\chi_\kappa$ is the dimension of representation of the symmetric group given by
$$ \chi_{\kappa}= k! \frac{\prod_{i<j}^{m}(k_i - k_j - i + j)}{\prod_{i=1}^{m}(k_i + m - i)!}.$$

The linearization coefficients, $g^\delta_{\kappa,\tau}(\beta)$, appear in the distribution for multivariate hypothesis-testing problems, see Sugiyama~\cite{sugiyama1970joint}, Pillai and Sugiyama~\cite{pillai1969non}, Ratnarajah et al. \cite{ratnarajah2004eigenvalues} and Shimizu and Hashiguchi~\cite{shimizu2022algorithm, shimizu2022numerical}.
An explicit representation for the coefficients $g^\delta_{\kappa, \tau}(1)$ with $\kappa$ of length one was proposed by Kushner \cite{kushner1988linearization}.
However, an algorithm for computing the general coefficients $g^\delta_{\kappa,\tau}(\beta)$ was unknown until recently.
Shimizu and Hashiguchi \cite{shimizu2022algorithm} provided an algorithm based on the expansion in terms of the elementary symmetric functions of Jack polynomials.
A different approach based on the expansion of power-sum symmetric functions was also developed by Hillier and Kan \cite{hillier2022properties}.
If $\kappa=(3,2)$, $\tau=(2,1)$, and $m=2$, the product of complex zonal polynomials is represented by 
\begin{align}
\label{g-coeff}
C_{(3,2)}(X)\,C_{(2,1)}(X)=\frac{5}{14}C_{(5,3)}(X)+\frac{5}{7}C_{(4,4)}(X).
\end{align}
Typically, a significant computation time is consumed for the complex zonal polynomials for large degrees. 
Efficient computation for obtaining $g^\delta_{\kappa,\tau}(\beta)$
is required.
In combinatorics, product expansions of symmetric polynomials, such as Macdonald and Schur polynomials, as in \eqref{g-coeff}, have been extensively studied.
The Schur polynomials $ \mathcal{S}_\kappa(X) $ can be defined as normalizing the complex zonal polynomials $ \mathcal{C}_\kappa(X) $ by
\begin{align}
\label{def-schur}
\mathcal{S}_\kappa(X) = \chi_\kappa^{-1} C_\kappa(X).
\end{align}
A single Schur polynomial has several determinantal representations, and those that allow high-precision computation in floating-point arithmetic have been discussed by Demmel and Koev \cite{demmel2006accurate} and Edelman et al. \cite{chan2019computing}.
The product of two Schur polynomials is represented by
\begin{align}
\label{prod-Schur}
S_\kappa(X) \, S_\tau(X) = \sum_{\delta \in P^{k+t}_m} c^\delta_{\kappa,\tau} S_\delta(X),
\end{align}
where $c^\delta_{\kappa,\tau}$ are determined by the Littlewood--Richardson rule, see Macdonald \cite{macdonald1998symmetric}.  
The coefficients, $c^\delta_{\kappa,\tau}$, are always non-negative integers.
If one of the polynomials in the product is indexed by a partition of length one (e.g., $( l(\kappa) = 1 $), then $c^\delta_{(k),\tau} = 1$. 
This is referred to as Pieri's rule. 
For real zonal polynomials, an explicit expression for $g^{\delta}_{{(k)}, \tau}(1)$ was provided by Kushner \cite{kushner1988linearization}.
Shimizu and Hashiguchi \cite{shimizu2022numerical} used Kushner's formula to compute the distribution of a real noncentral Roy's largest root statistic under a linear alternative.
From \eqref{def-schur} and \eqref{prod-Schur}, we obtain the Pieri formula for complex zonal polynomials.
\begin{lemma}
Let $k$ be a nonnegative integer and let $\delta=(\delta_1,\dots, \delta_l)$ and $\tau=(\tau_1,\dots,\tau_l)$ be partitions of $k+t$ and $t$, respectively.
Then we have 
 \begin{align}
 \label{g-coeffcient}
  C_{(k)}(X)\, C_\tau(X)
  =\sum_{\delta\in P^{k+t}_m} \frac{\chi_\tau}{\chi_\delta}C_\delta(X).
\end{align}
\end{lemma}
Similar to the Pieri formula, Corollary 4.2 in Naqvi \cite{naqvi2016product} provided a necessary and sufficient condition for $c^\delta_{\kappa,\tau}=1$.
If $\kappa \in P^k_2$ and $\tau \in P^t_2$, the coefficient $c^\delta_{\kappa,\tau}$ in (5) equals $1$ and it holds that
 \begin{align}
  \label{g-coeffcient-length2}
  C_{\kappa}(X)\, C_\tau(X)
  =\sum_{\delta\in P^{k+t}_2}  \frac{\chi_\kappa \chi_\tau}{\chi_\delta}C_\delta(X).
\end{align}
From $\chi_{(3,2)} = 5$, $\chi_{(2,1)} = 2$, and $\chi_{(5,3)} = 28$, it follows that 
$\chi_{(3,2)} \chi_{(2,1)} / \chi_{(5,3)}$ coincides with 
$g^{(5,3)}_{(3,2),(2,1)}(2) = 5/14$ in \eqref{g-coeff}.
The reduction formula of the computation for $c^\delta_{\kappa,\tau}$ was provided by Cho and Moon \cite{cho2011reduction}, and a special case of Richard Stanley's conjecture was proved by Naqvi \cite{naqvi2016product}. 
Some of these combinatorial results are useful for computing statistical distributions.

\section{Exact distribution of the noncentral complex Roy's largest root statistics}
\label{sec:03}
Let $\vec{x}_i$, $i = 1, \ldots, n$, be mutually independent distributed to the complex multivariate Gaussian distribution $CN_m(\vec{\mu}_i, \Sigma)$, where $\Sigma$ is an $m \times m$ positive-definite covariance matrix.  
Define the $m \times n$ matrix $X = [\vec{x}_1, \ldots, \vec{x}_n]$ and the mean matrix $M = [\vec{\mu}_1, \ldots, \vec{\mu}_n]$.  
The random matrix $W = XX^H$ is said to follow a complex noncentral Wishart distribution, denoted by $CW_m(n, \Sigma,\Omega)$, where the superscript $(\cdot)^H$ denotes the Hermitian transpose, and $\Omega = \Sigma^{-1} M M^H$.  
The central case is denoted by $CW_m(n, \Sigma)$.
The spectral decomposition of $W$ with $m>n$ is $W = E_1 \Lambda E_1^H$, 
where $\Lambda$ is an $n \times n$ diagonal matrix and the  $m \times n$ matrix $E_1$ is satisfied with $E_1^H E_1=I_n$.
The set of all  such $m\times n$ matrices $E_1$ with orthonormal columns is called the complex Stiefel manifold $CV_{n, m}$ defined by $CV_{n, m}=\{E_1 \in \mathbb{C}^{m \times n} \mid E_1^H E_1=I_n \}$,
where $n \le m$. 
If $m=n$, $U(m)=CV_{n, m}$ is the unitary group.

The Jacobian transformations for deriving the density of a singular Wishart matrix were given by Uling \cite{uhlig1994singular}.
He also presented a conjecture concerning the Jacobian transformations to derive the singular matrix variate beta yype I and II distributions, which was later proved by D\'iaz-Garc\'ia and Guti\'errez-J\'aimez \cite{diaz1997proof}.
These results are for the real case, whereas the corresponding results for the complex case were provided by D\'iaz-Garc\'ia and Guti\'errez-J\'aimez \cite{diaz2013distributions}.
In the same manner of D\'iaz-Garc\'ia et al. \cite{diaz1997proof}, the density of the noncentral complex singular Wishart matrix is given as 
\begin{align*}
f(W)=\frac{\pi^{n(n-m)}}{|\Sigma|^n \Gamma_n(n)}|\Lambda|^{n-m}\mathrm{etr}( -\Sigma^{-1} W) 
\mathrm{etr} ( -\Sigma^{-1} MM^H  ) {}_0F_1(n;\Omega \Sigma^{-1}W), 
\end{align*}
where $\mathrm{etr}(\cdot)=\exp(\mathrm{tr}(\cdot))$ and ${}_aF_b$ is the hyper geometric function of  a matrix of one argument (see Muirhead) \cite{muirhead1982aspects}.

Using the noncentral Wishart matrix, we define the matrix variate noncentral complex beta type I and II distributions. 
Let $B=(A_1+A_2)^{-1/2}A_1(A_1+A_2)^{-1/2}$ and $F=A_2^{-1/2}A_1A_2^{-1/2}$, where $A_1\sim CW_m(n,I_m,\Omega)$, $A_2\sim CW_m(p,I_m)$, and $p\geq m$. 
These matrices are called noncentral complex singular Beta and $F$ matrices.
The definitions of these matrices differ in the literature because the corresponding densities can be explicitly expressed. 
However, all cases for the eigenvalue distributions can be derived and are thus identical.
The eigenvalues of $B$ are denote by $1>\ell_1>\dots >\ell_{n_\mathrm{min}}>0$, $n_\mathrm{min}=\mathrm{min}\{n, m\}$. 
The eigenvalues of $F$ are given as $\ell_i/(1-\ell_i)$ ,$(i=1,\dots,n_\mathrm{min})$. 
Because the eigenvalue distribution of one can be derived from that of the other, we consider the density function and eigenvalue distribution of the Beta matrix.
\begin{proposition}
\label{prop-dist-B}
Let $A_1$ and $A_2$ be independent, where $A$ is $CW_m(n, I_m,\Omega)$ and $B$ is $CW_m(p,I_m)$ with $p \geq m> n$.
Then the density of $B=(A_1+A_2)^{-1/2}A_1(A_1+A_2)^{-1/2}$ is given as 
\begin{align}
\label{dist-B}
\nonumber
f(B)&=\frac{\pi^{n(n-m)}}{\Gamma_n(n)\Gamma_m(p)}|L|^{n-m}|I_m-B|^{p-m}\mathrm{etr} \left( - \Omega  \right)\\
&\times \int_{C>0}|C|^{n+p-m}\mathrm{etr}(-C)  {}_0F_1(n;\Omega C^{1/2}BC^{1/2})\, (dC),
\end{align}
where $B = E_1 L E_1^H$, $E_1 \in CV_{n,m}$, and $L$ is an $n \times n$ diagonal matrix.
The notation $C>0$ indicates that $C$ is Hermitian positive definite. 
\begin{proof}
The density functions of $A_1$ and $A_2$ are given as
\begin{align*}
f(A_1)&=\frac{\pi^{n(n-m)}}{\Gamma_{n}(n)}|\Lambda|^{n-m}\mathrm{etr} \left( - A_1 \right) 
\mathrm{etr} \left( - \Omega  \right)  {}_0F_1(n;\Omega A_1)\\
f(A_2)&=\frac{1}{\Gamma_m(p)}|A_2|^{p-m}
\mathrm{etr}(-A_2),
\end{align*}
respectively. 
Thus, we have 
\begin{align*}
f(A_1,A_2)=\frac{\pi^{n(n-m)}}{\Gamma_n(n)\Gamma_m(p)}|\Lambda|^{n-m}|A_2|^{p-m}
\mathrm{etr}(-(A_1+A_2))\mathrm{etr} \left( - \Omega  \right)  {}_0F_1(n;\Omega A_1).
\end{align*}
Using the transformation $A_1+A_2=C$ with Jacobian transformation of Lemma 3 in D\'iaz-Garc\'ia and Guti\'errez-J\'aimez \cite{diaz2013distributions}, we obtain the joint density of $B$ and $C$ as
\begin{align*}
f(B,C)=\frac{\pi^{n(n-m)}}{\Gamma_n(n)\Gamma_m(p)}|L|^{n-m}|C|^{n+p-m}
\mathrm{etr}(-C)|I_m-B|^{p-m}\mathrm{etr} \left( - \Omega  \right)  {}_0F_1(n;\Omega C^{1/2}B C^{1/2}).
\end{align*}
Integrating the above density with respect to $C > 0$, we obtain the desired result. 
\end{proof}
\end{proposition}
\begin{remark}
In the real and nonsingular case, Gupta and Nagar \cite{gupta2018matrix} also provided the final expression of the noncentral Beta matrix in the form of an integral representation. 
Diaz-Garcia et al. \cite{diaz2007noncentral} evaluated such integrals using the concept of ''symmetrised density" introduced by Greenacre \cite{greenacre1973symmetrlsed}.
 In Proposition~\ref{prop-dist-B}, a similar approach can be applied to evaluate the integral; however, because we are interested in the eigenvalue distribution, we do not pursue discussion on the densities of Beta and $F$ matrices. 
 \end{remark}
In the nonsingular case, the joint density of eigenvalues of a noncentral complex $F$-matrix was provided in James \cite{james1964distributions}, where the density is obtained with respect to the Lebesgue measure.
From Proposition~\ref{prop-dist-B}, we have the joint density of eigenvalues for the singular case. 
\begin{proposition}
\label{eigen-denisityB}
Let $A_1$ and $A_2$ be independent, where $A_1$ is $CW_m(n, I_m,\Omega)$ and $A_2$ is $CW_m(p,I_m)$ with $p \geq m> n$.
The joint density of eigenvalues $\ell_1,\dots,\ell_n$ of $B=(A_1+A_2)^{-1/2}A_1(A_1+A_2)^{-1/2}$ is given as 
\begin{align}
\label{dist-joint-eigenvalue}
\nonumber
f(\ell_1,\dots,\ell_n)=\ &\frac{\pi^{n(n-1)}\Gamma_m(n+p)}{\Gamma_n(n)\Gamma_m(p)\Gamma_n(m)}|L|^{m-n}|I_n-L|^{p-m}\prod_{i<j}^{n}(\ell_i-\ell_j)^2 \\
&\times \mathrm{etr} \left( - \Omega  \right) 
 \sum_{k=0}^{\infty}\sum_{\kappa \in P_n^k}\frac{(n+p)_\kappa C_\kappa(\Omega)C_\kappa(L)}{(n)_\kappa k! C_\kappa(I_m)}.
\end{align}
\begin{proof}
The Jacobian transformation of $B=E_1LE_1^H$ provided by Diaz-Garcia and Gutierrez-Sanchez \cite{diaz2013distributions} is given as 
\begin{align*}
(dB)=2^{-n}\pi^{-n}\prod_{i=1}^{n}|L|^{2(m-n)}\prod_{i<j}^{n}(\ell_i-\ell_j)^{2}(dL)\wedge(E_1^{H}dE_1). 
\end{align*}
Using $(dE_1)=\Gamma_n(m)/2^n\pi^{mn}(E^H_1dE_1)$ and Theorem 1 in Shimizu and Hashiguchi \cite{shimizu2021heterogeneous},
the function ${}_0F_1$ in \eqref{dist-B} is expanded as 
\begin{align}
\label{zonal-expansiton}
\int_{E_1\in CV_{n,m}}{}_0F_1(n;\Omega C^{1/2}E_1LE_1^H C^{1/2})\, (dE_1)=\sum_{k=0}^{\infty}\sum_{\kappa \in P_n^k}\frac{ C_\kappa(C\Omega)C_\kappa(L)}{(n)_\kappa k! C_\kappa(I_m)},
\end{align}
where $C=A_1+A_2$.
Then we have
\begin{align*}
\frac{\pi^{n(n-1)}}{\Gamma_n(n)\Gamma_m(p)\Gamma_n(m)}|L|^{m-n}|I_n-L|^{p-m}\prod_{i<j}^{n}(\ell_i-\ell_j)^{2}\mathrm{etr} \left( - \Omega  \right)\\
\times  \sum_{k=0}^{\infty}\sum_{\kappa \in P_n^k}\frac{ C_\kappa(L)}{(n)_\kappa k! C_\kappa(I_m)}\int_{C>0}|C|^{n+p-m}\mathrm{etr}(-C) C_\kappa(C\Omega)\, (dC).
\end{align*}
Integrating the above density with respect to $C>0$ using (2.7) in Ratnarajah et al. \cite{ratnarajah2004eigenvalues}, we obtain the desired results.
\end{proof}
\end{proposition}
The following lemma, first provided in the real case by Sugiyama \cite{sugiyama1967distribution}, is a useful integral formula for obtaining the distribution of the largest eigenvalue from the joint density of the eigenvalues.
Its extension to complex cases was given by Shimizu and Hashiguchi \cite{shimizu2021heterogeneous}.
   \begin{lemma}
   \label{Sugiyama-formula}
  Let $X_1=\mathrm{diag}(1,x_2,\dots, x_n)$ and $X_2=\mathrm{diag}(x_2,\dots, x_n)$ with $x_2>\cdots >x_n>0$; then the following equation holds. 
  \begin{align*}
\int_{1>x_2>\cdots >x_n>0}|X_2|^{a-n}C_\kappa(X_1)\prod_{i=2}^{n}(1-x_i)^2 \prod_{i<j}(x_i-x_j)^2 \prod_{i=2}^{n}dx_i\\
=(na+k)(\Gamma_n(n)/\pi^{n^2-n})\frac{\Gamma_n(a,\kappa)\Gamma_n(n)C_\kappa(I_n)}{\Gamma_n(a+n,\kappa)},
\end{align*}
where $\mathrm{\Re}(a)>n-1$.
\end{lemma}
From Proposition~\ref{eigen-denisityB} and Lemma~\ref{Sugiyama-formula}, we obtain the exact distribution of the noncentral complex Roy's largest root statistic.
\begin{theorem}
\label{prob-ell1}
Let $A_1$ and $A_2$ be independent, where $A_1$ is $CW_m(n, I_m,\Omega)$ and $A_2$ is $CW_m(p,I_m)$ with $p \geq m> n$.
Then the distribution function of $\ell_1$ of $B=(A_1+A_2)^{-1/2}A_1(A_1+A_2)^{-1/2}$ is given as 
\begin{align*}
\mathrm{Pr} (\ell_1<x)=\ &\frac{\Gamma_m(n+p)\Gamma_{n}(n)}{\Gamma_{m}(p)\Gamma_{n}(m+{n})}~\mathrm{etr}(-\Omega)\sum_{k=0}^{\infty}\sum_{\kappa\in P^k_{n}}\frac{(n+p)_\kappa}{(n)_\kappa}\frac{C_\kappa(\Omega)}{C_\kappa(I_m)k!}\\
& \times \sum_{t=0}^{\infty}\sum_{\tau\in P^t_{n}}\sum_{\delta \in P^{k+t}_{n}}g^\delta_{\kappa,\tau}(2) \frac{(m-p)_\tau }{t!} \frac{(m)_\delta \mathcal{C}_\delta (I_n)} {(m+n)_\delta}x^{mn+k+t}.
\end{align*}
\begin{proof}
This proof is similar to that of  Pillai and Sugiyama~\cite{pillai1969non}.
From \eqref{prod-zonal}, we write
\begin{align*}
|I_{n}-L|^{p-m} C_\kappa(L)&={_1F_0}(m-p,L)C_\kappa(L)\\
&=\sum_{t=0}^{\infty}\sum_{\tau\in P^t_{n}} \frac{(m-p)_\tau C_\tau (L)}{t!}C_\kappa(L)\\
&=\sum_{t=0}^{\infty}\sum_{\tau\in P^t_{n}}\sum_{\delta \in P^{k+t}_{n}}\frac{(m-p)_\tau}{t!}g^\delta_{\kappa,\tau}(2)C_\delta(L). 
\end{align*}
Therefore, the density function \eqref{dist-joint-eigenvalue} is represented by 
\begin{align*}
f(\ell_1,\dots, \ell_n)=&~C_0~\mathrm{etr}(-\Omega)|L|^{m-n}
\sum_{k=0}^{\infty}\sum_{\kappa\in P^k_{n}}\frac{(n+p)_\kappa}{n_\kappa}\frac{C_\kappa(\Omega)}{C_\kappa(I_m)k!}\\
&\times \prod_{i<j}^{n}(\ell_i-\ell_j)^2\sum_{t=0}^{\infty}\sum_{\tau\in P^t_n}\sum_{\delta \in P^{k+t}_n}g^\delta_{\kappa,\tau}(2) \frac{(m-p)_\tau C_\delta(L)}{t!},
\end{align*}
where 
\[
C_0=\frac{\pi^{n(n-1)}\Gamma_m(n+p)}{\Gamma_n(n)\Gamma_m(p)\Gamma_n(m)}.
\]
Translating $x_i=\ell_i/\ell_1$, $i=2,\dots,n$ and using Lemma~\ref{Sugiyama-formula}, we have
\begin{align*}
f(\ell_1)=&~C_0~\mathrm{etr}(-\Omega)\sum_{k=0}^{\infty}\sum_{\kappa\in P^k_{n}}\frac{(n+p)_\kappa}{n_\kappa}\frac{C_\kappa(\Omega)}{C_\kappa(I_m)k!}\sum_{t=0}^{\infty}\sum_{\tau\in P^t_{n}}\sum_{\delta \in P^{k+t}_{n}}g^\delta_{\kappa,\tau} (2)\frac{(m-p)_\tau }{t!}\\
&\times \ell_1^{mn+k+t-1}  \int_{1>x_2>\cdots>x_{n}>0}|X_2|^{m-n}\prod_{i=2}^{n}(1-x_i)\prod_{2\leq i<j}^{n}(\ell_i-\ell_j)^2C_\delta(X_1)\prod_{i=2}^{n}dx_i\\
=&~\frac{\Gamma_m(n+p)\Gamma_{n}(n)}{\Gamma_{m}(p)\Gamma_{n}(m+{p})}~\mathrm{etr}(-\Omega)\sum_{k=0}^{\infty}\sum_{\kappa\in P^k_{n}}\frac{(n+p)_\kappa}{(n)_\kappa}\frac{C_\kappa(\Omega)}{C_\kappa(I_m)k!}\\
& \times \sum_{t=0}^{\infty}\sum_{\tau\in P^t_{n}}\sum_{\delta \in P^{k+t}_{n}}g^\delta_{\kappa,\tau} (2)\frac{(m-p)_\tau }{t!}\ell_1^{mn+k+t-1} (mn+k+t)\frac{(m)_\delta \mathcal{C}_\delta (I_n)}
{(m+n)_\delta},
\end{align*}
where matrices $X_1$ and $X_2$ are defined in Lemma~ \ref{Sugiyama-formula}.
Finally, by integrating $f(\ell_1)$ with respect to $\ell_1$, we obtain the desired result. 
\end{proof}
\end{theorem}
The joint density of the eigenvalues in the nonsingular case was derived by James \cite{james1964distributions}. 
In a manner similar to Theorem 1, the distribution of the largest eigenvalue can be obtained from the joint density of the eigenvalues.
The following expression generalizes the results to cover both nonsingular and singular cases.
\begin{corollary}
\label{corollary-prob-rank1}
Let $A_1$ and $A_2$ be independent, where $A_1$ is $CW_m(\Omega, n,I_m)$ and $A_2$ is $CW_m(p,I_m)$ with $p \geq m$, $\mathrm{rank}(\Omega)=1$ and $\theta_1$ is the largest eigenvalue of $\Omega$, then we have
\begin{align}
\label{prob-rank1}
\nonumber
\mathrm{Pr} (\ell_1<x)&=\frac{\Gamma_m(n+p)\Gamma_{n_\mathrm{min}}(n_\mathrm{min})}{\Gamma_{m}(p)\Gamma_{n_\mathrm{min}}(m+n)}~\mathrm{etr}(-\Omega)\sum_{k=0}^{\infty}\frac{(n+p)_\kappa \theta^k_1}{(n)_\kappa(m)_\kappa}\\
& \times \sum_{t=0}^{\infty}\sum_{\tau\in P^t_{n_\mathrm{min}}}\sum_{\delta \in P^{k+t}_{n_\mathrm{min}}}\frac{\chi_\tau}{\chi_\delta} \frac{(m-p)_\tau }{t!} \frac{(n_\mathrm{max})_\delta \mathcal{C}_\delta (I_{n_\mathrm{min}})} {(m+n)_\delta}x^{mn+k+t},
\end{align}
where $n_\mathrm{max}=\mathrm{max}\{n,m\}$ and $n_\mathrm{min}=\mathrm{min}\{n,m\}$.
\begin{proof}
From $\mathrm{rank}(\Omega)=1$, we have $C_\kappa(\Omega) = C_{(k)}(\Omega) = \theta_1^k$.  
Moreover, according to (3.29) of Drensky et al. \cite{drensky2014computing}, for partitions $\kappa$ with only one part, $C_{(k)}(I_m)$ is given by $C_{(k)}(I_m) = (m)_k/k!$.
\end{proof}
\end{corollary}
If $A_1 $ holds rank two in Theorem \ref{prob-ell1}, then from \eqref{g-coeffcient-length2}, the distribution of the largest eigenvalue $\ell_1$ is represented as
\begin{align*}
\mathrm{Pr} (\ell_1<x)&=\frac{\Gamma_m(n+p)\Gamma_{n_\mathrm{min}}(n_\mathrm{min})}{\Gamma_{m}(p)\Gamma_{n_\mathrm{min}}(m+n)}~\mathrm{etr}(-\Omega)\sum_{k=0}^{\infty}\sum_{\kappa\in P^k_{2}}\frac{(n+p)_\kappa}{(n)_\kappa}\frac{C_\kappa(\Omega)}{C_\kappa(I_m)k!}\\
& \times \sum_{t=0}^{\infty}\sum_{\tau\in P^t_2}\sum_{\delta \in P^{k+t}_2}\frac{\chi_\kappa\chi_\tau}{\chi_\delta} \frac{(m-p)_\tau }{t!} \frac{(n_\mathrm{max})_\delta \mathcal{C}_\delta (I_2)} {(m+n)_\delta}x^{mn+k+t}. 
\end{align*}
The above distribution, under the generalized alternative, does not require the algorithm for obtaining the linearization coefficients proposed by Shimizu and Hashiguchi \cite{shimizu2022algorithm}. 
Compared to Theorem \ref{prob-ell1}, it can be computed more efficiently using the explicit formula for Schur polynomials.
We also obtain the hth moment of $\ell_1$.
\begin{corollary}
\label{corollary: prob-rank1-moment}
Under the same condition of Corollary~\ref{corollary-prob-rank1}, the $h$th moment of $\ell_1$ is given as 
\begin{align}
\label{prob-rank1-moment}
\nonumber
\mathrm{E}[\ell^h_1] &=\frac{\Gamma_m(n+p)\Gamma_{n_\mathrm{min}}(n_\mathrm{min})}{\Gamma_{m}(p)\Gamma_{n_\mathrm{min}}(m+n)}~\mathrm{etr}(-\Omega)\sum_{k=0}^{\infty}\frac{(n+p)_\kappa \theta^k_1}{(n)_\kappa(m)_\kappa}\\
& \times \sum_{t=0}^{\infty}\sum_{\tau\in P^t_{n_\mathrm{min}}}\sum_{\delta \in P^{k+t}_{n_\mathrm{min}}}\frac{\chi_\tau}{\chi_\delta} \frac{(m-p)_\tau }{t!} \frac{(n_\mathrm{max})_\delta \mathcal{C}_\delta (I_{n_\mathrm{min}})} {(m+n)_\delta}\frac{mn+k+t}{mn+k+t+h}.
\end{align}
\end{corollary}

\section{Numerical computation}
\label{sec:04}
\subsection{Computation of the distribution}
In this section, we calculate the distribution of a noncentral complex Roy's largest root statistic.
The results in the previous section are presented as infinite series.  
By noting the properties of the Pochhammer symbol, the summation with respect to $t$ in \eqref{prob-rank1} becomes a finite sum up to $n(p - m)$.
Then the truncated distribution up to the $K$th degree is given by
\begin{align}
\label{trunc-dist}
\nonumber
F_K(x)&=\frac{\Gamma_m(n+p)\Gamma_{n_\mathrm{min}}(n_\mathrm{min})}{\Gamma_{m}(p)\Gamma_{n_\mathrm{min}}(n+m)}~\mathrm{etr}(-\Omega)\sum_{k=0}^{K}\frac{(n+p)_k \theta^k_1}{(n)_k(m)_k}\\
& \times \sum_{t=0}^{n(p-m)}\sum_{\tau\in P^t_{n_\mathrm{min}}}\sum_{\delta \in P^{k+t}_{n_\mathrm{min}}}\frac{\chi_\tau}{\chi_\delta} \frac{(m-p)_\tau }{t!} \frac{(n_\mathrm{max})_\delta \mathcal{C}_\delta (I_{n_\mathrm{min}})} {(m+n)_\delta}x^{mn+k+t}.
\end{align}
Figure~\ref{fig1} illustrates the graph of the distribution in \eqref{trunc-dist}, with $\theta_1 = 9$, $n = 2$, $m=4$ and $p = 12$, for $K = 5, 10, 20$. 
Numerical computations were performed using Mathematica $14.2.1$ on a computer~(MacBook Pro macOS Sonoma, ver. 14.4.1, M2 Chip with 24GB RAM).
As $K$ increases, the cumulative probability of the truncated distribution approaches $1$.  
The other derived results can also be computed by applying the same truncation approach.
The $95$th percentile of the largest eigenvalue obtained from $10^6$ Monte Carlo simulations is $0.782$. 
Substituting this into the exact distribution yields an upper tail probability of $0.95$, confirming that a truncation at $K = 20$ is sufficient.

\begin{figure}[H]
\label{fig1}
\begin{center}
\includegraphics[width=8cm]{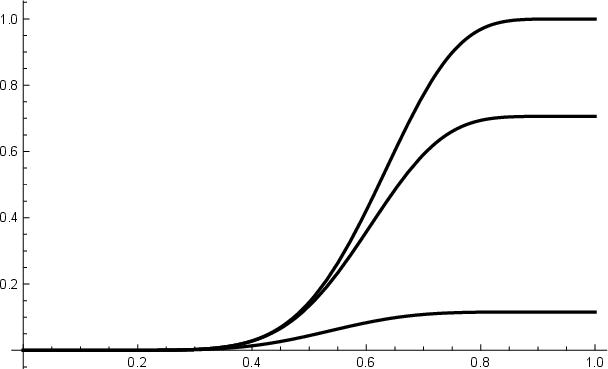}
\rlap{\raisebox{32.5ex}{\kern-23.em{\large $F_K(x)$}}}%
\rlap{\raisebox{30ex}{\kern-4em{$K=20$}}}%
\rlap{\raisebox{22ex}{\kern-4em{$K=10$}}}%
\rlap{\raisebox{6ex}{\kern-4em{$K=5$}}}%
\rlap{\raisebox{.15cm}{\kern0cm{\large $x$}}}
\caption{Truncated distribution with $\theta_1=9$, $n=2$, $m=4$, and $p=12$.}
 \label {fig1}
\end{center}
\end{figure}
Table~\ref{moment} shows the mean, variance, skewness, and kurtosis of $\ell_1$ with $n=2$, $p=18$ and $\theta_1=2$, computed using the result of Corollary~\ref{corollary: prob-rank1-moment}.
The mean and kurtosis increase while the variance and skewness become larger as $m$ increases. 
\begin{table}[H]
  \begin{center}
\caption{Expectation and variance of $\ell_1$ with $n=2$, $p=18$, and $\theta_1=2$} 
\label{moment}
{\begin{tabular}
{@{}ccccc@{}} \toprule
$m$& 
$\mathrm{Mean}$&
$\mathrm{Variance}$&
$\mathrm{Skewness}$&
$\mathrm{kurtosis}$\\  \toprule
     5   & 0.402 & 0.00978 &  0.177 & 2.82 \\
    10   & 0.649 & 0.00789 & -0.182 & 2.85 \\
    15   & 0.867 & 0.00349 & -0.688 & 3.50 \\
\noalign{\smallskip}\hline
\end{tabular}}
 \end{center}
\end{table}
\subsection{Computation of the test power}
In this subsection, we discuss the test of equality for complex mean vectors from $d$ groups.  
Let $\vec{x}_{ij}$ be distributed as the $i$-th population $CN_m(\vec{\mu}_i, \Sigma)$, $i=1,\dots, d$, $j=1,\dots, n_i$, where $N-d\geq m$ and $N=\sum_{i=1}^{d}n_i$.
We consider testing the equality of the mean vectors as follows:
\begin{align}
\label{test-signal}
H_0\!:  \vec{\mu}_1=\cdots =\vec{\mu}_d \ \text{vs.} \ H_1\!:\vec{\mu} \neq 0 \ \text{for some} j .
\end{align}
The between-group and within-group covariance matrices are defined by
$H=\sum_{i=1}^{d}n_i(\bar{\vec{x}}_i-\bar{\vec{x}})(\bar{\vec{x}}_i-\bar{\vec{x}})^\top, E=\sum_{i=1}^{d}\sum_{j=1}^{n_i}(\vec{x}_{ij}-\bar{\vec{x}}_i)(\vec{x}_{ij}-\bar{\vec{x}}_i)^\top$, respectively, where $\bar{\vec{x}}_i=n^{-1}_i\sum_{j=1}^{n_i}\vec{x}_{ij}$ and $\bar{\vec{x}}=\sum_{i=1}^{d}\sum_{j=1}^{n_i}\vec{x}_{ij}/n$.
The matrices $H$ and $E$ are distributed as $CW_m(n_H,\Sigma,\Omega)$ and $CW_m(n_E, \Sigma)$, respectively, where $n_H=d-1$ and $n_E=N-d$.
Here, $\Omega$ is defined as $\Omega = \Sigma^{-1} \sum_{i=1}^{d} n_i (\vec{\mu}_i - \vec{\bar{\mu}})(\vec{\mu}_i - \vec{\bar{\mu}})^\top$, where $\vec{\bar{\mu}} = N^{-1} \sum_{i=1}^{d} \vec{\mu}_i$.
To select a more appropriate test statistic for complex MANOVA, it is important to compare the powers of several test statistics discussed in the introduction.
We compute the power of the complex Roy's test statistic using the theoretical results.
If $\Omega=O$ in \eqref{prob-rank1}, the null distribution required to obtain the critical value is expressed by 
\begin{align*}
\mathrm{Pr} (\ell_1<x)=\frac{\Gamma_m(n_H+n_E)\Gamma_{n_H}(n_\mathrm{min})}{\Gamma_{m}(n_E)\Gamma_{n_\mathrm{min}}(m+n_H)} x^{mn_H}
{}_2F_1(m-n_E, m;m+n_H;xI_{n_H}).
\end{align*}
The null distribution does not require the computation of linearization coefficients, making it easier to compute than the noncentral case. 
Table \ref{table:power comparison} shows the exact and approximate powers of Roy's test at the nominal significance level of $0.05$. 
\begin{table}[H]
\caption{Exact and approximate power for the linear case}
\begin{center}
\begin{tabular}{@{}ccc ccc ccc@{}} \toprule
$m$ & $n_i$ & $d$ & 
\multicolumn{2}{c}{$\theta_1 = 10$} & 
\multicolumn{2}{c}{$\theta_1 = 20$} & 
\multicolumn{2}{c}{$\theta_1 = 30$} \\ 
\cmidrule(lr){4-5} \cmidrule(lr){6-7} \cmidrule(lr){8-9}
& & & Exact & Approx. & Exact & Approx. & Exact & Approx. \\ \midrule
3  & 3 & 3 & 0.430 & 0.435 & 0.770 & 0.808 &  0.926 &   0.950    \\ 
3  & 5 & 3 &    0.677   &   0.672   &    0.961 &   0.967    &    0.997   & 0.998      \\
3  & 3 & 5 &   0.407    &   0.338  &  0.783  &   0.766    &    0.946   &   0.947    \\
7  & 4 & 4 &   0.208    &  0.107   &   0.452    &  0.387  &    0.678   &   0.667    \\
10 & 7 & 3 & 0.271 & 0.184 & 0.590 & 0.579 &  0.823  &  0.854     \\ 
\bottomrule
\label{table:power comparison}
\end{tabular}
\end{center}
\end{table}
The approximate power was obtained using the approximate distribution proposed by Dharmawansa \cite{dharmawansa2019roy}. 
We can confirm that the approximation performs well in low dimensions, but its accuracy decreases as the dimension increases.

\section*{Conclusion}
In this study, we discussed the linear coefficients $g^\delta_{\kappa, \tau}(2)$ that arise when expressing the product of complex zonal polynomials.
By applying the Littlewood--Richardson rule and Pieri's formula, which are well known in combinatorics, we provided an explicit expression for the linear coefficients.
We derived the distribution of the noncentral complex Roy's statistic using Pieri's formula when the noncentrality matrix has rank one and performed a numerical evaluation.

In combinatorics, algorithms and explicit expressions for computing Littlewood--Richardson coefficients associated with Schur polynomials, and their generalizations have been developed 
\cite{cho2011reduction,naqvi2016product}.
These developments have a strong potential for use in problems involving the distributional computations of test statistics, particularly those explored in this study.
Combinatorial results may be applied to the numerical evaluation of the distribution under general alternatives in Theorem~\ref{prob-ell1}, which remains a topic for future work.
\section*{Acknowledgments}
This work was supported by JSPS KAKENHI (Grant Numbers 25K17300 and 23K11016).

\bibliographystyle{elsarticle-harv}
\bibliography{ref}




\end{document}